\documentclass[a4paper,12pt]{article}

\PassOptionsToPackage{hyphens}{url}
\usepackage{graphicx}
\usepackage{amssymb}
\usepackage[fleqn]{amsmath}
\usepackage{booktabs}
\usepackage[T1]{fontenc} 
\usepackage{lmodern}
\usepackage{hyperref}
\usepackage{subfig}
\usepackage{natbib}
\usepackage{bm}
\usepackage{enumitem}
\usepackage{tikz}

\renewcommand{\vec}[1]{\boldsymbol{#1}}
\newcommand{\field}[1]{\mathbb{#1}}  
\newcommand{\R}{\field{R}} 
\newcommand{\Uset}{\mathcal{U}} 
\newcommand{\transp}{^{^{\top}}}
\newcommand{\onenorm}[1]{\left|\left|#1\right|\right|_1}
\newcommand{\twonorm}[1]{\left|\left|#1\right|\right|_2}
\newcommand{\infnorm}[1]{\left|\left|#1\right|\right|_\infty}

\newtheorem{theorem}{Theorem}[section]
\newenvironment{proof}[1][Proof]{\noindent\textbf{#1.} }{\ \rule{0.5em}{0.5em}}

\hypersetup{
    colorlinks=true,
    citecolor=black,
    urlcolor=blue,
    linkcolor=black,
    pdfborder={0 0 0},
}

\begin{document}

\title{Robust fractional programming}

\author{Bram L. Gorissen \\ \\
\textit{\small Tilburg University, Department of Econometrics and Operations Research} \\
\textit{\small 5000 LE Tilburg, Netherlands} \\
\textit{\small {\tt b.l.gorissen@tilburguniversity.edu}}}
\date{}
\maketitle

\begin{abstract}
We extend Robust Optimization to fractional programming, where both the objective and the constraints contain uncertain parameters. Earlier work did not consider uncertainty in both the objective and the constraints, or did not use Robust Optimization. Our contribution is threefold. First, we provide conditions to guarantee that either a globally optimal solution, or a sequence converging to the globally optimal solution, can be found by solving one or more convex optimization problems. Second, we identify two cases for which an exact solution can be obtained by solving a single optimization problem: (1) when uncertainty in the numerator is independent from the uncertainty in the denominator, and (2) when the denominator does not contain an optimization variable. Third, we show that the general problem can be solved with an (iterative) root finding method. The results are demonstrated on a return-on-investment maximization problem, data envelopment analysis, and mean-variance optimization. We find that the robust optimal solution is only slightly more robust than the nominal solution. As a side-result, we use Robust Optimization to show that two existing methods for solving fractional programs are dual to each other.
\end{abstract}

\begin{tikzpicture}[remember picture,overlay]
\node[anchor=south,yshift=5pt] at (current page.south) {\fbox{\parbox{\dimexpr\textwidth-\fboxsep-\fboxrule\relax}{\footnotesize This is an author-created, un-copyedited version of an article published in Journal of Optimization Theory and Applications \href{http://dx.doi.org/10.1007/s10957-014-0633-4}{DOI:10.1007/s10957-014-0633-4}.}}};
\end{tikzpicture}

\section{Introduction}

A fractional program (FP) is an optimization problem, where the objective is a fraction of two functions. It can be used for an economical trade-off such as maximizing return/investment, maximizing return/risk or minimizing cost/time \citep{Schaible1983}. A comprehensive overview of FP papers, containing over 550 references which include many applications, is given by \citet{Schaible1982}. More up-to-date references can be found in \citep{Stancu2013}, which also refers to six preceding bibliographies by the same author.

More often than not, the parameters in an optimization problem are affected by uncertainty. Robust Optimization (RO) is about solving optimization problems with uncertain data in a computationally tractable way (see, e.g., \citet{BenTal,Bertsimas2011}). The key concept is that a solution has to be feasible for all realizations of the uncertain parameters, which are assumed to reside in convex uncertainty regions.

Sometimes, the objective is the maximum of finitely many fractions, and the feasible region is a convex set. This is called a generalized FP (see, e.g., \cite{Barros1996,crouzeix1991} for solution methods). A generalized FP with infinitely many fractions in the objective was solved by \cite{Lin2005} using a cutting plane method, that uses a set of finitely many fractions that is extended in each step. They do not mention that their method can be used to deal with uncertain data, and do not use existing results from RO. Our work can be seen as an alternative approach, where we also deal with uncertainty in the constraints.

The Lagrange dual of a robust FP was studied by \cite{Jeyakumar2013}, extending a result by \cite{Beck2009}. They assume that the uncertainty in the numerator of the objective is independent from the uncertainty in the denominator. The dual is tractable when the numerator, denominator and constraints are linear, and the uncertainty regions are ellipsoidal or finite sets of scenarios. In this paper, we focus on the primal problem. Nevertheless, in Section \ref{sec:singlesolverrun}, we obtain and extend the list of tractable duals.

The aim of this paper is to combine FP and RO, to provide a comprehensive overview of the solution methods, and to investigate the improvement of RO on numerical examples. First, we provide conditions that guarantee that a globally optimal solution, or a sequence that converges to the globally optimal solution can be found by solving one or more convex problems. The importance of these conditions is illustrated with a numerical example from literature. Second, we identify two cases for which an exact solution can be obtained by solving a single optimization problem. Third, we show that the general problem can be solved with an (iterative) root finding method.

In Section \ref{sec:nonrobustfp}, we outline two existing solution methods for FPs, and present a new result showing that the two approaches are each others duals. In Section \ref{sec:ro}, we present existing results in RO, that will be used for FPs as well. Our main contribution is given in Section \ref{sec:solving-rfp}. The results are demonstrated on a return-on-investment maximization problem, data envelopment analysis, and mean-variance optimization in Section \ref{sec:numex}.
\section{Solving Nonrobust Fractional Programs}\label{sec:nonrobustfp}
In this section, we present two existing methods to solve FPs, and show that these methods are dual to each other. To the best of our knowledge, this duality result is new.
Consider the following general formulation of an FP:
\begin{align*}
\mbox{(FP)} \quad
\min_{\vec{x} \in \R^n} \quad \frac{f(\vec{x})}{g(\vec{x})} \quad
\mbox{s.t.} \quad  h_i(\vec{x}) \leq 0, \; \forall i \in I.
\end{align*}
We will assume that the constraint index set $I$ is finite, that $f$ is convex and non-negative and that $g$ is concave and positive over the feasible region. When the functions $f$, $g$ and $h_i$ are affine, (FP) is a {\it linear} fractional program:
\begin{align*}
\mbox{(LFP)}       \quad
\min_{\vec{x} \in \R^n}  \quad \frac{b_0 + \vec{b} \transp \vec{x}}{c_0 + \vec{c} \transp \vec{x}} \quad \mbox{s.t.} \quad d_{0i} + \vec{d_i} \transp \vec{x} \leq 0, \quad \; i \in I.
\end{align*}
Charnes and Cooper show that (LFP) can be reformulated as an (equivalent) LP, by making the substitutions $\vec{y} = \vec{x} /(c_0 + \vec{c} \transp \vec{x})$ and $t=1 /(c_0 + \vec{c} \transp \vec{x})$ \citep{CharnesCooper}:
\begin{align*}
\mbox{(CC-LFP)} \quad \min_{t \in \R_+,\vec{y} \in \R^n} \quad b_0 t + \vec{b} \transp \vec{y} \quad \mbox{s.t.} \quad d_{0i} t + \vec{d_i} \transp \vec{y} \leq 0, \; \forall i \in I, \quad  c_0 t + \vec{c} \transp \vec{y} = 1.
\end{align*}
An optimal solution of (LFP) is obtained from an optimal solution of (CC-LFP) by computing $\vec{x}=\vec{y}/t$. \cite{Schaible1974} shows that a similar substitution ($\vec{y} = \vec{x} / g(\vec{x})$, $t = 1/g(\vec{x})$) transforms (FP) into an equivalent convex programming problem:
\begin{align}
\mbox{(Schaible-FP)} \quad \min_{t \in \R_{++},\vec{y} \in \R^n} \quad t f\left(\frac{\vec{y}}{t}\right) \quad \mbox{s.t.} \quad t g\left(\frac{\vec{y}}{t}\right)   \geq 1, \quad t h_i\left(\frac{\vec{y}}{t}\right) \leq 0, \; \forall i \in I. \notag
\end{align}
This is indeed a convex problem, since the perspective function $p(\vec{y},t) := t f\left(\vec{y}/t\right)$ is jointly convex on $\R^n \times \R_{+}$ when $f$ is convex on $\R^n$. Furthermore, an optimal $\vec{x}$ is obtained from $\vec{x}=\vec{y}/t$. Schaible also shows that, if the constraint $t g(\vec{y}/t) \geq 1$ is formulated as an equality constraint: $t g(\vec{y}/t) = 1$, it is not necessary to require $f$ to be positive \citep{Schaible1974}.

A different solution approach uses the auxiliary parameterized optimization problem $F(\alpha)$, defined as $\min_{\vec{x} \in \R^n} \{ f(\vec{x}) - \alpha g(\vec{x}) : h_i(\vec{x}) \leq 0 \, \; \forall i \in I \}$. The objective value of (FP) is at least $\alpha$ if and only if $F(\alpha) \geq 0$ \citep{Dinkelbach1967}. So, the objective of (FP) equals the largest $\alpha$ such that $F(\alpha) \geq 0$:
\begin{align*}
\mbox{(P-FP)}   \quad \max_{\alpha \in \R_+} \{  \alpha : \min_{\vec{x} \in X} \{ f(\vec{x}) - \alpha g(\vec{x}) \} \geq 0 \},
\end{align*}
where $X$ denotes the feasible region $\{ \vec{x} \in \R^n : h_i(\vec{x}) \leq 0, \; \forall i \in I \}$. The usual way of solving this problem is by finding the root of $F$, since the corresponding $\vec{x}$ is optimal for (FP) \citep{Dinkelbach1967}. This is usually done with a Newton-like algorithm, where there is some freedom in choosing the next iteration point \citep{ChenFP2009}. The root of $F$ is unique, since $F$ is monotonically decreasing in $\alpha$. The parameteric program $F(\alpha)$ is convex when $g$ is affine or when $f$ is non-negative on the feasible region (since then, only non-negative values for $\alpha$ needs to be considered). For these cases, the Newton method to find the root of $F$ was described by \cite{Dinkelbach1967}, which creates a monotonically decreasing sequence, that converges superlinearly and often (local) quadratically to a root of $F$ \citep{Schaible1976}.

We now show that these approaches are dual to each other when $f$ is non-negative. The proof for affine $g$ and possibly negative $f$ is similar.
\begin{theorem}\label{thm:dual}\textup{
Assume that $f$ is convex and non-negative on $X$, $g$ is concave on $X$, $X$ is closed and convex, and the optimal value of (FP) is attained. Then (Schaible-FP) and (P-FP) are dual to each other, and strong duality holds.}
\end{theorem}
\begin{proof}
First note that the following reformulation is equivalent to (P-FP):
\begin{align*}
\mbox{(RP-FP)}   \quad \max_{\alpha \in \R_+} \{  \alpha :  f(\vec{x}) - \alpha g(\vec{x}) \geq 0, \; \forall \vec{x} \in X \}.
\end{align*}
The remainder of the proof is based on the theory ``primal worst is dual best'' introduced by \cite{Beck2009}. They assume that $X$ is compact and convex. Since we have assumed $X$ to be closed and convex and since an optimal solution of (FP) is attained, compactness can be achieved by intersecting $X$ with a box that includes the optimal solution. Additionally, Beck and Ben-Tal assume that the constraint in (RP-FP) is convex in $\alpha$ and concave in $\vec{x}$, which indeed holds. 
For fixed $\vec{x}$, (RP-FP) is an LP with the following dual:
\begin{align*}
\mbox{(D-LP)}     \quad \min_{t \in \R_+} \{  t f(\vec{x}) :  t g(\vec{x}) \geq 1 \}.
\end{align*}
While (RP-FP) is robust since the constraint has to hold for all $\vec{x}$ in $X$, the constraint in the optimistic counterpart of (D-LP) has to hold for a single $\vec{x}$:
\begin{align*}
\mbox{(OD-LP)}    \quad \min_{t \in \R_+, \vec{x} \in X} \{  t f(\vec{x}) :  t g(\vec{x}) \geq 1 \}.
\end{align*}
(RP-FP) and (OD-LP) are dual to each other \citep{Beck2009}. Strong duality holds, since $(\vec{x},t)$ is a Slater point for (OD-LP) for sufficiently large $t$. It is obvious that (OD-LP) and (Schaible-FP) are equivalent, since $t=0$ is infeasible for (OD-LP).
\end{proof}

\section{Robust Optimization}\label{sec:ro}
There are currently two generic methods that deal with an infinite number of constraints. The first method is applicable to both linear and nonlinear constraints, while the second method can only be applied to robust LPs.

The first (``constraint wise'') approach uses conic duality \citep{BenTal} or Fenchel duality \citep{BenTal2011}. The vector $\vec{x}$ in $\R^n$ satisfies the following infinite number of constraints:
\begin{align*}
&&h_i(\vec{a_i}, \vec{x}) \leq 0, \;\forall \vec{a_i} \in \R^L : \tau_{ij}(\vec{a_i}) \leq 0, \; \forall j \in J,
\end{align*}
if and only if there exists $u_{ij} \in \R_+$ and $\vec{v_{ij}} \in \R^L$ ($j \in J$, $J$ being a finite set), such that $\vec{x}$ satisfies the following convex constraint:
\begin{align}
&&\sum_{j \in J} u_{ij} \tau_{ij}^*\left( \frac{ \vec{v_{ij}} }{ u_{ij} }\right) - (h_i)_*\left(\sum_{j \in J} \vec{v_{ij}}, \vec{x}\right) \leq 0, \label{eq:expl:fenchel}
\end{align}
where $\tau_{ij}^*(\vec{s}) = \sup_{\vec{a_i} \in \R^L}\{ \vec{s} \transp \vec{a_i} - \tau_{ij}(\vec{a_i})\}$ and $(h_i)_*(\vec{s}, \vec{x}) = \inf_{\vec{a_i} \in \R^L}\{ \vec{s} \transp \vec{a_i} - h_i(\vec{a_i},\vec{x})\}$ are the convex and concave conjugates of $\tau_{ij}$ and $h_i$, respectively. This approach requires the constraint to be concave in $\vec{a_i}$, the functions $\tau_{ij}$ to be convex, and ri(dom $h_i(\cdot,\vec{x})$) $\cap$ ri($\Uset_i$) $\neq \emptyset$ for all $\vec{x} \in \R^n$, where $\Uset_i := \{\vec{a_i} \in \R^L : \tau_{ij}(\vec{a_i}) \leq 0, \; \forall j \in J \}$. This approach yields a tractable formulation for many constraints and many uncertainty sets (see Tables 1 and 2 in \citep{BenTal2011}), even if the conjugates do not have closed-form expressions. To give an impression of the broad applicability of this method, let us cite some examples, for which it provides a tractable reformulation. For uncertainty sets, one could have a norm-bounded (e.g., box or ball), polyhedral or conic representable set, or a generic set defined by (convex) power functions, exponential functions, negative logarithms, or any function for which the convex conjugate exists. Constraints could be linear or quadratic in the uncertain parameter. There are some operations on functions that preserve the availability of a tractable expression of the conjugate. One is multiplication with a non-negative scalar: the concave conjugate of $t f(\vec{a_0},\vec{y}/t)$ (with respect to $\vec{a_0}$) for $t\geq0$, is the perspective of the concave conjugate of $f$: $t f_*(\vec{s} / t, \vec{y} / t)$. Another one is when $h_i$ is the sum of two functions: $h_i = h_{i1} + h_{i2}$. Suppose closed-form conjugates exist for $h_{i1}$ and $h_{i2}$ separately; then:
\begin{align}
&& (h_i)_*(\vec{s},\vec{x}) = \max_{\vec{s_1} \in \R^L,\vec{s_2} \in \R^L} \left\{ (h_{i1})_*(\vec{s_1},\vec{x}) + (h_{i2})_*(\vec{s_2},\vec{x}) : \vec{s_1} + \vec{s_2} = \vec{s} \right\}. \label{eq:fenchel:sum}
\end{align}
When substituting \eqref{eq:fenchel:sum} into \eqref{eq:expl:fenchel}, the max operator in \eqref{eq:fenchel:sum} may be omitted, since if the resulting constraint holds for some $\vec{s_1}$ and $\vec{s_2}$, then it surely holds for the maximum. This example shows that a closed-form expression for the conjugate is indeed not always required.

The second method solves any robust LP with a convex uncertainty region via its Lagrange dual \citep{Gorissen2013convunc}. The transformation from the primal to the dual is a three step procedure. First, the dual of the nonrobust LP is formulated, where the uncertain parameters are assumed to be known. Second, instead of enforcing the constraints for all realizations of the uncertain parameters (``robust counterpart''), the constraints of the dual have to hold for a single realization of the uncertain parameters (``optimistic counterpart''). So, the uncertain parameters are added to the set of optimization variables. The last step is to reformulate the nonconvex optimistic counterpart to an equivalent convex optimization problem. The optimal solution of the resulting problem can be translated to an optimal solution of the original robust LP via the KKT vector.

In the remainder, we provide reformulations based on \eqref{eq:expl:fenchel}, but the reader should be aware that the other approach may be useful when all functions involved are linear.

\section{Solving Robust Fractional Programs}\label{sec:solving-rfp}
In this section, we show how to solve (R-FP). It is our aim to obtain Robust Counterparts (RCs), that can be solved with existing Robust Optimization methods. First, we formulate conditions that give raise to convex optimization problems. Under these conditions, a globally optimal solution can be found by solving a single convex optimization problem (Sections \ref{sec:singlesolverrun} and \ref{sec:singlesolverrun2}). These conditions also guarantee that, in the general case (Section \ref{sec:generalcase}), a root finding method produces a sequence of convex optimization problems, whose solutions converge to a globally optimal solution. The results of this section are summarized in Tables \ref{tbl:results1}-\ref{tbl:results3}.

\subsection{Robust Formulation and Assumptions}
The uncertain parameters, denoted by $\vec{a_i}$, are assumed to lie in sets $\Uset_i \subset \R^L$, which we define using functions $\tau_{ij} : \R^L \to \R$:
\begin{align}
&& \Uset_i := \{\vec{a_i} \in \R^L : \tau_{ij}(\vec{a_i}) \leq 0, \; \forall j \in J \}, \label{eq:uset}
\end{align}
where $J$ is a finite set. In the RC of (FP), the constraints have to be satisfied by all realizations of the uncertain parameters:
\begin{align}
\mbox{(R-FP)} \quad \min_{\alpha \in \R,\vec{x} \in \R^n} \quad ( \alpha ) \quad
\mbox{s.t.} \quad &\frac{f(\vec{a_0}, \vec{x})}{g(\vec{a_0}, \vec{x})} - \alpha \leq 0, \; \forall \vec{a_0} \in \Uset_0 \notag \\
& h_i(\vec{a_i}, \vec{x}) \leq 0, \; \forall \vec{a_i} \in \Uset_i, \; \forall i \in I. \label{eq:feasub}
\end{align}
Note that the uncertainty is specified constraint-wise, which is possible even if the parameters in different constraints are correlated \citep[p. 12]{BenTal}.
We make the following assumptions:
\begin{enumerate}[label=(\alph*)]
\item $\tau_{ij}$ are convex and the sets $\Uset_i$ are convex and compact,
\item $f$ and $h_i$ are convex in $\vec{x}$ for every fixed value of $\vec{a_i}$ in $\Uset_i$,
\item $g$ is concave in $\vec{x}$ for every fixed value of $\vec{a_i}$ in $\Uset_i$,
\item $f$ and $h_i$ are concave in $\vec{a_i}$ for every feasible $\vec{x}$,
\item $g$ is convex in $\vec{a_i}$ for every feasible $\vec{x}$,
\item $g(\vec{a_0},\vec{x}) > 0$ for every $\vec{a_0}$ in $\Uset_0$ and every feasible $\vec{x}$, and
\item $f(\vec{a_0}, \vec{x}) \geq 0$ for at least one $\vec{a_0} \in \Uset_0$ and for every feasible $\vec{x}$.
\end{enumerate}
The last assumption is not necessary if $g$ is biaffine, i.e., when $g$ is affine in each parameter when the other parameter is fixed, of which we show the consequences in Section \ref{subsec:biaffine}. In robust {\it linear} programming, the assumption that $\Uset_i$ is compact and convex is made without any loss of generality \citep[p.~12]{BenTal}. For robust FP, compactness is not a restriction since the functions $h_i$ are continuous and the constraints do not contain strict inequalities. So, the problem remains unchanged if the uncertainty region is replaced with its closure. However, requiring $\Uset_i$ to be convex is a restriction (unless $f$, $g$ and $h$ are affine in $\vec{a_i}$), that is necessary for using existing results in RO.

Assumptions (d) and (e) are made solely because they are required by generic methods to derive a tractable RC. There are some examples where the RC can be derived even though these conditions are not fulfilled, e.g., if the uncertainty region is the convex hull of a limited number of points and the constraint is convex in the uncertain parameter, for a conic quadratic program with implementation error, or when the $\mathcal{S}$-lemma or a sums of squares result can be applied \citep{BenTal,BenTal2011,BenTal2014hidden,BertsimasPolynomials}. In these cases, assumptions (d) and (e) are not necessary.

In literature, a problem is solved, that does not satisfy these requirements \citep{Lin2005}. While the authors claim that this did not affect their computations, and that they find the global optimum, in Appendix \ref{sec:lin2005} we show that their solution is suboptimal.

\subsection{Special Case: Uncertainty in the Numerator is Independent of the Uncertainty in the Denominator}\label{sec:singlesolverrun}
Suppose that the uncertainty in the numerator of the objective is decoupled from the uncertainty in the denominator:
\begin{align*}
\mbox{(R-S1)} \quad \min_{\alpha \in \R_+,\vec{x} \in \R^n} \quad (\alpha) \quad \mbox{s.t.} \quad &\frac{f(\vec{a_0}, \vec{x})}{g(\vec{a_0}', \vec{x})} - \alpha \leq 0, \; \forall \vec{a_0} \in \Uset_{0}, \; \forall \vec{a_0}' \in \Uset_{0}' \\
& h_i(\vec{a_i}, \vec{x}) \leq 0, \; \forall \vec{a_i} \in \Uset_i, \; \forall i \in I.
\end{align*}
We claim that (R-S1) is equivalent to the RC of the Schaible reformulation:
\begin{align}
\mbox{(R-Schaible-S1)} \quad \min_{\alpha \in \R_+,t \in \R_{++},\vec{y} \in \R^n}\quad (\alpha) \quad \mbox{s.t.} \quad & t f\left(\vec{a_0}, \frac{\vec{y}}{t}\right) - \alpha \leq 0, \; \forall \vec{a_0}  \in \Uset_{0} \notag \\
& t g\left(\vec{a_0}', \frac{\vec{y}}{t}\right)   \geq 1, \; \forall \vec{a_0}' \in \Uset_{0}' \label{s-s1-eq} \\
& t h_i\left(\vec{a_i}, \frac{\vec{y}}{t}\right) \leq 0, \; \forall \vec{a_i}  \in \Uset_i, \; \forall i \in I. \notag
\end{align}
Clearly, an optimal solution of (R-Schaible-S1) exists, for which $t = \sup_{\vec{a_0}' \in \Uset_{0}'} 1/g(\vec{a_0}',\vec{y}/t)$. Equivalence between (R-S1) and (R-Schaible-S1) readily follows from the substition $\vec{x}=\vec{y}/t$, that converts a feasible solution of one problem to a feasible solution of the other problem.

This result extends \citep{Jeyakumar2013}. They provide the dual of (R-S1), and show that strong duality holds. In case $f$, $g$ and $h_i$ are linear, they show that the dual of (R-S1) is tractable, when the uncertainty region is ellipsoidal, or consists of a finite set of scenarios. The resulting problems can also be obtained from our work, by applying the solution method by \cite{Gorissen2013convunc} to (R-Schaible-S1). In addition to ellipsoids or scenarios, our method works with any convex uncertainty region, such as a polyhedral set or a conic quadratic representable set. A similar result was found by \cite{Kaul1986}, but that result is wrong (see Appendix \ref{sec:kaul86}).

We provide a reformulation of (R-Schaible-S1) if the RO method using conjugates is used (eq. \eqref{eq:expl:fenchel}). The resulting equivalent problem becomes:
\begin{align*}
\min \quad (\alpha) \quad \mbox{s.t.} \quad & \sum_{j \in J} u_{0j} \tau_{0j}^*\left(\frac{\vec{v_{0j}} }{ u_{0j}}\right) - t f_*\left(\frac{\sum_{j \in J} \vec{v_{0j}}}{t}, \frac{\vec{y}}{t}\right) - \alpha \leq 0 \\
& \sum_{j \in J} u_{0'j} \tau_{0'j}^*\left(\frac{\vec{v_{0'j}} }{ u_{0'j}}\right) - t g_*\left(\frac{\sum_{j \in J} \vec{v_{0'j}}}{t}, \frac{\vec{y}}{t}\right)   \geq 1  \\
&\sum_{j \in J} u_{ij} \tau_{ij}^*\left(\frac{\vec{v_{ij}} }{ u_{ij}}\right) - t (h_i)_*\left(\frac{\sum_{j \in J} \vec{v_{ij}}}{t}, \frac{\vec{y}}{t}\right) \leq 0, \; \forall i \in I \\
& \alpha \in \R_+,t \in \R_{++},\vec{u}\in \R^{(|I|+2)\times|J|},\vec{v}\in \R^{(|I|+2)\times|J|\times L},\vec{y} \in \R^n.
\end{align*}

\subsection{Special Case: the Denominator Does Not Depend on the Optimization Variable \texorpdfstring{$\vec{x}$}{x}}\label{sec:singlesolverrun2}
If the optimization variables do not appear in the denominator, (R-FP) is equivalent to \citep[cf.][Ex.~30]{BenTal2011}:
\begin{align*}
\mbox{(R-S2)} \quad
\min_{\alpha \in \R_+,\vec{x} \in \R^n} \quad (\alpha) \quad \mbox{s.t.} \quad & f(\vec{a_0}, \vec{x}) - \alpha g(\vec{a_0}) \leq 0, \; \forall \vec{a_0} \in \Uset_0 \\
& h_i(\vec{a_i}, \vec{x}) \leq 0, \; \forall \vec{a_i} \in \Uset_i, \; \forall i \in I.
\end{align*}
Note that $g$ indeed does not depend on $\vec{x}$. (R-S2) can be solved via the following equivalent convex reformulation using \eqref{eq:expl:fenchel}, and standard techniques for the conjugate of the sum of two functions:
\begin{align*}
\min \quad (\alpha) \quad \mbox{s.t.} \quad & \sum_{j \in J} t_{0j} \tau_{0j}^*\left(\frac{\vec{v_{0j}} }{ t_{0j} } \right) - f_*\left(\vec{s} + \sum_{j \in J} \vec{v_{0j}}, \vec{x}\right) + \alpha g^*\left(\frac{\vec{s} }{ \alpha }\right) \leq 0 \\
& \sum_{j \in J} t_{ij} \tau_{ij}^*\left(\frac{\vec{v_{ij}} }{ t_{ij} } \right) - (h_i)_*\left(\vec{\sum_{j \in J} v_{ij}}, \vec{x}\right) \leq 0, \; \forall i \in I \\
& \alpha \in \R_+,\vec{t}\in \R^{(|I|+1)\times|J|}_+,\vec{s} \in \R^L,\vec{v} \in \R^{(|I|+1)\times|J|\times L},\vec{x} \in \R^n.
\end{align*}

\subsection{General Case}\label{sec:generalcase}
We now show how to solve the general problem (R-FP) using the following parametric problem:
\begin{align}
F(\alpha) := \min_{\vec{x},w}  \quad (w)  \quad \mbox{s.t.} \quad & f(\vec{a_0}, \vec{x}) - \alpha g(\vec{a_0}, \vec{x}) \leq w, \; \forall \vec{a_0} \in \Uset_0 \label{eq:compf} \\
& h_i(\vec{a_i}, \vec{x}) \leq 0, \; \forall \vec{a_i} \in \Uset_i, \; \forall i \in I, \notag
\end{align}
which is a convex optimization problem, since we only have to solve it for $\alpha \in \R_+$. Let $\alpha^*$ be a root of $F$. Lin and Sheu show that an optimal solution of (R-FP) is the minimizer $\vec{x}$ of $F(\alpha^*)$, if the feasible region for $\vec{x}$ is compact \citep{Lin2005}. We assume from now on that the constraint functions $h_i$ define a compact feasible region. Moreover, they show that $F(\alpha) < 0$ if and only if $\alpha > \alpha^*$. Lin and Sheu do not use results from RO to arrive at the deterministic reformulation \eqref{eq:Fexpl}. Instead, they replace the set $\Uset_0$ with a finite set, approximate $F(\alpha)$ with an entropic regularization method, and iteratively generate a sequence $\tilde{\alpha_k}$ that converges to $\alpha^*$. The approximation becomes more accurate as the root of $F$ is approached. The reason why they approximate $F(\alpha)$ is because they claim that computing its value is difficult. Our approach is to solve $F(\alpha)$ using RO, which inherently produces tractable problems. The following convex reformulation using \eqref{eq:expl:fenchel} is equivalent:
\begin{align}
F(\alpha) = \min \quad (w)  \quad \mbox{s.t.} \quad & \sum_{j \in J} t_{0j} \tau_{0j}^*\left(\frac{\vec{v_{0j}} }{ t_{0j} }\right) - f_*\left(\vec{s} + \sum_{j \in J} \vec{v_{0j}}, \vec{x}\right) + \alpha g^*\left(\frac{\vec{s} }{ \alpha}, \vec{x}\right) \leq w \label{eq:Fexpl} \\
& \sum_{j \in J} t_{ij} \tau_{ij}^*\left(\frac{\vec{v_{ij}} }{ t_{ij} }\right) - (h_i)_*\left(\sum_{j \in J} \vec{v_{ij}}, \vec{x}\right) \leq 0, \; \forall i \in I \notag \\
& \vec{t} \in \R^{(|I|+1) \times |J|}_+,\vec{s} \in \R^L,\vec{v} \in \R^{(|I|+1) \times |J| \times L},w \in \R,\vec{x} \in \R^n. \notag
\end{align}
Since $F$ is monotonically decreasing in $\alpha$, as for FPs and generalized FPs, existing root-finding methods can be used. We mention a few of these that produce a sequence $\{ \alpha_k \}$ which converges to $\alpha^*$:
\begin{enumerate}[label=(\alph*)]
\item The bisection method. Bounds on the interval that contain $\alpha^*$ are: \begin{align}
&& & \alpha_{LB} := \min_{\vec{x} \in \R^n}  \{ f(\vec{\bar{a}_0},\vec{x}) / g(\vec{\bar{a}_0},\vec{x}) : h_i(\vec{a_i}, \vec{x}) \leq 0, \; \forall \vec{a_i} \in \Uset_i, \; \forall i \in I \} \label{eq:lb} \\
&& & \alpha_{UB} := \sup_{\vec{a_0} \in \Uset_0} f(\vec{a_0},\vec{x}) / g(\vec{a_0},\vec{x}), \label{eq:ub}
\end{align}
where \eqref{eq:lb} is computed for a fixed $\vec{\bar{a}_0}$ from $\Uset_0$, and \eqref{eq:ub} is computed for some  $\vec{x}$ that is (robust) feasible, i.e., for an $\vec{x}$ that satisfies \eqref{eq:feasub}. These bounds can be computed relatively easily using the Schaible reformulation. If the lower bound \eqref{eq:lb} is hard to compute due to the ``for all'' quantifier in the constraints, it may be computed for fixed $\vec{\bar{a}_i}$ from $\Uset_i$. Since $F(\alpha_{LB}) \geq 0$ and $F(\alpha_{UB}) \leq 0$, and since $F$ is clearly nonincreasing, $\alpha^*$ lies in $[\alpha_{LB}, \alpha_{UB}]$. The middle point of this interval is $\alpha_k := 0.5(\alpha_{LB} + \alpha_{UB})$. By evaluating $F(\alpha_k)$, the width of the interval that contains $\alpha^*$ can be halved: if $F(\alpha_k) > 0$, then set $\alpha_{LB} = \alpha_k$, otherwise set $\alpha_{UB} = \alpha_k$. By increasing $k$ by 1 and repeating this procedure, a series $\{ \alpha_k \}$ is constructed, that converges to $\alpha^*$.
\item The Dinkelbach type algorithm by \cite{crouzeix1985}, adjusted for infinitely many ratios. The method starts with $k=0$ and $\alpha_k = \sup_{\vec{a_0} \in \Uset_0} f(\vec{a_0},\vec{x}) / g(\vec{a_0},\vec{x})$ for some feasible $\vec{x}$. Then $F(\alpha_k)$ is computed, with maximizer $\vec{x_k}$. If $F(\alpha_k) < 0$, then the next $\alpha$ is determined by $\alpha_{k+1} := \max_{\vec{a_0} \in \Uset_0} f(\vec{a_0}, \vec{x_k}) / g(\vec{a_0}, \vec{x_k})$. Computing $\alpha_{k+1}$ requires solving an FP. The method proceeds by increasing $k$ by 1, and again computing $F(\alpha_k)$. If the feasible region for $\vec{x}$ is compact, then the series $\{ \alpha_k \}$ converges linearly to $\alpha_k$.
\item The same as method (b), except that the right-hand side of \eqref{eq:compf} is multiplied with $g(\vec{a_0}, \vec{x_k})$: $f(\vec{a_0}, \vec{x}) - \alpha g(\vec{a_0}, \vec{x}) \leq w g(\vec{a_0}, \vec{x_k})$. This may increase the speed of convergence for the same complexity of computation \citep{crouzeix1986}.
\item The same as method (b) or (c), except that the $\vec{a_0}$ that maximizes $F(\alpha_k)$ is used to compute $\alpha_{k+1}$, instead of solving a new optimization problem. The worst case $\vec{a_0}$ in the computation of $F(\alpha_k)$ can be recovered without much computational effort. Thus, $\alpha_{k+1} = f(\vec{a_0}, \vec{x_k}) / g(\vec{a_0}, \vec{x_k})$ is computed more efficiently than in the method (b) or (c). Additional work is required to ensure convergence of $\{ \alpha_k \}$ to $\alpha^*$ \citep[Section~5]{crouzeix1991}.
\end{enumerate}
Let $\vec{x_k}$ be the maximizer of $F(\alpha_k)$. If a root finding method finds the root in a finite number of steps, then an exact solution of (R-FP) is found. Otherwise, Crouzeix et al. show that, if the sequence $\{ \alpha_k \}$ converges to $\alpha^*$, then any convergent subsequence of $\{ \vec{x_k} \}$ converges to the optimal solution $\vec{x^*}$ of (R-FP) \citep[Theorem~4.1c]{crouzeix1985}.

\subsection{Consequences when the Denominator is Biaffine}\label{subsec:biaffine}
The assumption that the numerator is positive, ensures that the objective value of (R-FP) is positive over the feasible region. Consequently, we could assume $\alpha \in \R_+$; this would produce a convex optimization problems. If $g$ is biaffine, then the resulting problems are also convex for $\alpha < 0$. We shall discuss the results to each of the three aforementioned cases separately. For the {\it first special case} (Section \ref{sec:singlesolverrun}), the restriction that the numerator is positive may be dropped only if the denominator does not contain an uncertain parameter. Then, (R-S1) and (R-Schaible-S1) are equivalent if $\alpha \in \R_+$ is replaced with $\alpha \in \R$ and \eqref{s-s1-eq} is stated as an equality (cf. \cite{Schaible1974}). The reason why the denominator may not contain an uncertain parameter, is because $t = 1/g(\vec{x},\vec{a_{01}})$ is not possible for multiple $\vec{a_{01}}$.

For the {\it second special case} (Section \ref{sec:singlesolverrun2}), the denominator only depends on $\vec{a_0}$, so ``biaffine'' in the title of this subsection should be read as ``affine''. When (R-S2) is solved for $\alpha \in \R$, the restriction that the numerator is positive, may be dropped. 

For the {\it general case} (Section \ref{sec:generalcase}), no changes need to be made to drop the restriction that the numerator is positive.

\begin{table}[h]
	\centering
	\captionsetup{width=0.8\textwidth}
  \caption{Tractable cases when uncertainty in the numerator is independent of the uncertainty in the denominator. $\ell$ denotes an affine function.}
  \label{tbl:results1}
	\begin{tabular}{lllll}
	\toprule
$f$ & $g$ & $h_i$ & sgn($f$) & (R-Schaible-S1) \\
   \midrule
$f(\vec{a_0},\vec{x})$     & $g(\vec{a_0}',\vec{x})$     & $h_i(\vec{a_i},\vec{x})$   & $\geq 0$ & no modifications                                 \\
$f(\vec{a_0},\vec{x})$     & $\ell(\vec{x})$       & $h_i(\vec{a_i},\vec{x})$   & any      & $\alpha \in \R$, \eqref{s-s1-eq} as an equality  \\
  \bottomrule
	\end{tabular}
\end{table}

\begin{table}[h]
	\centering
	\captionsetup{width=0.8\textwidth}
  \caption{Tractable cases when the denominator does not depend on $\vec{x}$. $\ell$ denotes an affine function.}
  \label{tbl:results2}
	\begin{tabular}{lllll}
	\toprule
$f$ & $g$ & $h_i$ & sgn($f$) & (R-S2) \\
   \midrule
$f(\vec{a_0},\vec{x})$     & $g(\vec{a_0})$    & $h_i(\vec{a_i},\vec{x})$   & $\geq 0$ & no modifications  \\
$f(\vec{a_0},\vec{x})$     & $\ell(\vec{a_0})$ & $h_i(\vec{a_i},\vec{x})$   & any      & $\alpha \in \R$   \\
  \bottomrule
	\end{tabular}
\end{table}

\begin{table}[h]
	\centering
	\captionsetup{width=0.8\textwidth}
  \caption{Tractable cases for the general case. $\ell$ denotes a biaffine function.}
  \label{tbl:results3}
	\begin{tabular}{llll}
	\toprule
$f$ & $g$ & $h_i$ & sgn($f$) \\
   \midrule
$f(\vec{a_0},\vec{x})$     & $g(\vec{a_0},\vec{x})$      & $h_i(\vec{a_i},\vec{x})$  & $\geq 0$  \\
$f(\vec{a_0},\vec{x})$     & $\ell(\vec{a_0},\vec{x})$   & $h_i(\vec{a_i},\vec{x})$  & any       \\
  \bottomrule
	\end{tabular}
\end{table}

\section{Numerical Examples}\label{sec:numex}
In this section, we test our method on three examples: a multi-item newsvendor problem (Section \ref{sec:mine}), mean-variance optimization (Section \ref{sec:mvo}), and data envelopment analysis (Section \ref{sec:dea}).

\subsection{Multi-item Newsvendor Example}\label{sec:mine}
In \cite{Gorissen2013convunc}, a multi-item newsvendor problem is solved by minimizing the investment cost under the condition that at least a certain expected profit is made.

We show how to directly optimize the expected return on investment for this example. Let us first recapitulate the problem. The newsvendor buys $Q_i$ units of item $i$ at the beginning of the day. Each item has its associated ordering cost $c_i$, selling price $v_i$, salvage price $r_i$, and unsatisfied demand loss $l_i$. We assume $r_i \leq v_i+l_i$. During the day the newsvendor faces a demand $d_i$, resulting in a profit of $\min \{ v_i Q_i + l_i(Q_i-d_i) - c_i Q_i, v_i d_i + r_i(Q_i-d_i) - c_i Q_i \}$. The demand is not known in advance, but there are finitely many demand scenarios $d_{is}$ ($s$ in $S$) that occur with (uncertain) probability $p_{is}$, independently of other items.

The problem of maximizing expected return on investment can be formulated as:
\begin{align}
\mbox{(R-NV)} \quad
\max_{\vec{Q} \in \R^{|I|}_+,\vec{u} \in \R^{|I| \times |S|}} \quad & \min_{\vec{p_i} \in \Uset_i} \quad \frac{ \sum_{i \in I}\sum_{s \in S} p_{is} u_{is}  }{  \sum_{i \in I} c_i Q_i   } \notag \\
\mbox{s.t.} \quad 
                                & u_{is} + \left( c_i - r_i \right) Q_i \leq d_{is} \left( v_i - r_i \right), \; \forall i \in I, \; \forall s \in S \notag \\
                                & u_{is} + \left( c_i - v_i - l_i \right) Q_i \leq - d_{is} l_i, \; \forall i \in I, \; \forall s \in S, \notag
\end{align}
where $u_{is}$ is the contribution to the profit of item $i$ in scenario $s$, and the convex and compact uncertainty regions $\Uset_i$ are defined using the Matusita distance, which is a $\phi$--divergence measure \citep{BTHWMR2011}:
\begin{align*} && \Uset_i = \left\{ \vec{p}_i \in \R^{|S|}_+ : \sum_{s \in S} p_{is} = 1, \sum_{s \in S} \left| \left(\hat{p}_{is}\right)^\alpha -  \left(p_{is}\right)^\alpha \right|^{1/\alpha} \leq \rho, \; \forall i \in I \right\}. \end{align*}
Note that the assumptions (a)-(e) and (g) are always fulfilled, and that assumption (e) is fulfilled if at least one item is bought. (R-NV) can be classified under the first special case (Section \ref{sec:singlesolverrun}). Since all functions are affine and the denominator is certain, the Schaible reformulation and the Charnes-Cooper reformulation are equivalent. (R-NV) is therefore equivalent to:
\begin{align}
\mbox{(R-CC-NV)} \quad \quad \quad \quad \notag \\
\max_{\vec{Q} \in \R^{|I|}_+,t \in \R_+,\vec{u} \in \R^{|I|\times|S|}} \quad  & \min_{\vec{p_i} \in \Uset_i} \quad  \sum_{i \in I}\sum_{s \in S} p_{is} u_{is} \notag \\
\mbox{s.t.} \quad 
                                & u_{is} + \left( c_i - r_i \right) Q_i \leq d_{is} \left( v_i - r_i \right) t, \; \forall i \in I, \; \forall s \in S \notag \\
                                & u_{is} + \left( c_i - v_i - l_i \right) Q_i \leq - d_{is} l_i t, \; \forall i \in I, \; \forall s \in S \notag \\
                                & \sum_{i \in I} c_i Q_i = 1, \notag
\end{align}
where $\vec{Q}$ and $\vec{u}$ in (R-CC-NV) have to be divided by $t$ to obtain the $\vec{Q}$ and $\vec{u}$ in (R-NV). (R-CC-NV) is a linear program with a convex uncertainty region that we solve via its dual, as outlined in the introduction. The last of the reformulation steps, a substitution, is not necessary since the uncertainty only appears in the objective. Let $x_{is}$, $y_{is}$ and $z$ be the dual variables of (CC-NV); then the optimistic dual (OD-CC-NV) is given by:
\begin{align}
\mbox{(OD-CC-NV)} \quad
\min \quad  & z \notag \\
\mbox{s.t.} \quad                    & x_{is} + y_{is} = p_{is}, \; \forall i \in I, \; \forall s \in S \notag \\
                                     & \sum_{s \in S} \{ (c_i - r_i)  x_{is} + (c_i-v_i-l_i) y_{is} \} + c_i z \geq 0, \; \forall i \in I \label{od-cc-nvQ} \\
                                     & \sum_{i \in I}\sum_{s \in S} - d_{is} (v_i-r_i) x_{is} + d_{is} l_i y_{is} \geq 0                                              \label{od-cc-nvT} \\
                                     & \sum_{s \in S} p_{is} = 1, \; \forall i \in I \notag \\
                                     & \sum_{s \in S} \left| \left(\hat{p}_{is}\right)^\alpha -  \left(p_{is}\right)^\alpha \right|^{1/\alpha} \leq \rho, \; \forall i \in I  \notag \\
                                     & \vec{p} \in \R^{|S|}_+, \vec{x} \in \R^{|I|\times|S|}_+, \vec{y} \in \R^{|I|\times|S|}_+,z\in\R. \notag
\end{align}
The optimal value of (OD-CC-NV) is the robust expected return on investment. The corresponding optimal order quantities $\vec{Q}$ can be derived from the KKT vector of (OD-CC-NV), by dividing its elements associated with \eqref{od-cc-nvQ} by the element corresponding to \eqref{od-cc-nvT}. This is the same as dividing $\vec{Q}$ by $t$ in (R-CC-NV) to undo the Charnes-Cooper transformation.

We solve the problem for the same data as \cite{Gorissen2013convunc} with AIMMS 3.11 (Paragon Decision Technology, the Netherlands) and KNITRO 7.0 (Zienna Optimization LLC, USA) with its default settings. Computation errors for negative $p_{is}$ were avoided by using $|p_{is}|^\alpha$ instead of $(p_{is})^\alpha$. We take $\rho = 0$ to obtain the nominal solution, while $\rho = 0.03$ for the robust solution. Solutions were obtained in less than 0.01 seconds. When the probabilities are as expected ($\rho = 0$), the expected return on investment of the nominal solution is 0.297, while for the robust solution it is 0.285. When $\rho = 0.03$ and the worst case probabilities occur for the nominal solution, i.e., the probabilities that minimize the expected return on investment for the nominal solution, the objective value drops to 0.211, while for the robust solution it drops to 0.214. So, the solution indeed becomes more robust, but the difference with the nominal solution is small. We verify if the decision maker could have done better, if he knew beforehand which probability vector realizes. This done by optimizing the nominal model ($\rho=0$), while setting the probabilty estimates $\hat{p}_{is}$ equal to the worst case probabilities for the robust solution. This gives the so-called perfect hindsight solution. The objective value is as low as 0.214. So, even though the robust objective could deteriorate substantially, there is no other solution that performs better.

\subsection{Mean-variance Optimization}\label{sec:mvo}
We are to present an example that involves a trade-off between mean and variance. This trade-off is commonly used in portfolio optimization, including the Modern Portfolio Theory (MPT) founded by \cite{Markowitz1952}, where the goal is to select the right mix of assets. In contrast to MPT, we do not impose that the expected returns on the assets and the covariance matrix are fully known. Instead, we assume that finitely many scenarios $s$ (in $S$) for the future can be identified along with unknown probabilities of occurence $p_s$, which are estimated by $\hat{p}_s$. The return of asset $i$ in scenario $s$ is a constant $r_{is}$, so when $x_i$ units of money are invested in asset $i$, the return in scenario $s$ is given by $u_s = \sum_{i \in I} r_{is} x_i$ (possibly negative). The expected return and variance are given by:
\begin{align}
	&&\mathbb{E}(\textrm{return})   &:= \sum_{s \in S} p_s u_s  \notag \\
	&&\mathrm{Var}(\textrm{return}) &:= \sum_{s \in S} p_s \left( u_s - \mathbb{E}( \textrm{return}) \right)^2 \label{varreturn1} \\
	&&                              &= \sum_{s \in S} p_s u_s^2 - \left(\sum_{s \in S} p_s u_s \right)^2. \label{varreturn2}
\end{align}
To remain in the minimization framework, the objective is to minimize the variance-to-mean ratio (or the dispersion index). The robust optimization problem is given by:
\begin{align*}
\mbox{(R-I)} \quad
\min_{\vec{x} \in \R^{|I|}_+,\vec{u} \in \R^{|S|}} \quad &\max_{\vec{p} \in \Uset}  \quad  \frac{\sum_{s \in S} p_s u_s^2 - \left(\sum_{s \in S} p_s u_s \right)^2}{\sum_{s \in S} p_s u_s}       \\
\mbox{s.t.} \quad	       & u_s = \sum_{i \in I} r_{is}  x_i, \; \forall s \in S \\
                           & \sum_{i \in I} x_i = C \\
                           & \sum_{s \in S} p_s u_s > 0, \; \forall \vec{p} \in \Uset.
\end{align*}
The last two constraints ensure that $C$ units of money are invested, and that this model has a feasible solution only if the expected profit is positive. The numerator is convex in $u_s$ (from \eqref{varreturn1}) and concave in $p_s$ (from \eqref{varreturn2}). The denominator is clearly concave in $p_s$ and convex in $u_s$.  Moreover, the numerator is non-negative and the denominator is positive on the feasible region. For the uncertainty region we use the modified $\chi^2$--distance as a $\phi$--divergence measure, which can be justified by statistical theory \citep{BTHWMR2011}:
\begin{align*} && \Uset_0 := \left\{ \vec{p}  \in \R^{|S|}_+ : \sum_{s \in S} p_s = 1, \sum_{s \in S} \frac{ \left( p_s - \hat{p}_s \right)^2 } { \hat{p}_s }  \leq \rho  \right\}. \end{align*}
(R-I) is not one of the special cases, so we solve this problem using the general method. In order to formulate \eqref{eq:Fexpl} explicitly, we first derive some conjugate functions. The conjugate for $f$ is from \citet[Ex.~25]{BenTal2011}.
\begin{align*}
& f(\vec{p},\vec{u}) := \sum_{s \in S} p_s u_s^2 - \left(\sum_{s \in S} p_s u_s \right)^2 && f_*(\vec{v},\vec{u}) = \sup_z \{ - \frac{z^2}{4} : u_s^2 + u_s z = v_s, \; \forall s \in S \} \\
& g(\vec{p},\vec{u}) := \sum_{s \in S} p_s u_s                                    && g^*(\vec{v},\vec{u}) = \begin{cases}  0,      &\mbox{if } v_s = u_s, \; \forall s \in S \\
                                                                                         \infty, & \mbox{otherwise} \end{cases} \\
& \tau_{i1}(\vec{p}) := \max_s \{ -p_s \}                           && \tau_{i1}^*(\vec{v}) =  \begin{cases}  0,      &\mbox{if } v_s \leq 0, \; \forall s \in S \mbox{ and } \sum_{s \in S} v_s \geq -1 \\
                                                                                         \infty, & \mbox{otherwise} \end{cases} \\
& \tau_{i2}(\vec{p}) := \sum_{s \in S} p_s - 1                              && \tau_{i2}^*(\vec{v}) =  \begin{cases}  1,      &\mbox{if } v_s =  1, \; \forall s \in S \\
                                                                                         \infty, & \mbox{otherwise} \end{cases} \\
& \tau_{i3}(\vec{p}) := 1 - \sum_{s \in S} p_s                              && \tau_{i3}^*(\vec{v}) =  \begin{cases} -1,      &\mbox{if } v_s = -1, \; \forall s \in S \\
                                                                                         \infty, & \mbox{otherwise} \end{cases} \\
& \tau_{i4}(\vec{p}) := \sum_{s \in S} \frac{ \left( p_s - \hat{p}_s \right)^2 } { \hat{p}_s }  - \rho  && \tau_{i4}^*(\vec{v}) = \rho + \sum_{s \in S} \hat{p}_s \left( \frac{1}{4} v_s^2 + v_s \right).
\end{align*}
Plugging in these formulas in \eqref{eq:Fexpl} yields $F(\alpha)=$
\begin{align}
                \min \quad w \quad \mbox{s.t.} \quad & t_{02} - t_{03} + \rho t_{04} + \sum_{s \in S} \hat{p}_s \left( \frac{(v_{04})_s^2}{4t_{04}} + (v_{04})_s \right) + \frac{z^2}{4} \leq w \label{eq:Fsubs1} \\
                & (v_{01})_s \leq 0 \quad (v_{11})_s \leq 0, \; \forall s \in S \notag \\
                & \sum_{s \in S} (v_{01})_s \geq -t_{01} \quad \sum_{s \in S} (v_{11})_s \geq -t_{11}  \notag \\
                & u_s^2 + u_s z = u_s \alpha + (v_{01})_s + t_{02} - t_{03} + (v_{04})_s, \; \forall s \in S \label{eq:Fsubs2} \\
                & t_{12} - t_{13} + \rho t_{14} + \sum_{s \in S} \hat{p}_s \left( \frac{(v_{14})_s^2}{4t_{14}} + (v_{14})_s \right) < 0 \notag \\
                & u_s + (v_{11})_s + t_{12} - t_{13} + (v_{14})_s = 0, \; \forall s \in S \notag \\
                & u_s = \sum_{i \in I} r_{is}  x_i,  \; \forall s \in S \notag \\
                & \sum_{i \in I} x_i = C \notag \\
                & \vec{t} \in \R^{2 \times 4}_+,\vec{u} \in \R^{|S|},\vec{v} \in \R^{2 \times 4},w \in \R,\vec{x} \in \R^{|I|},z \in \R. \notag
\end{align}
This problem is not convex because of the product $u_s z$ in \eqref{eq:Fsubs2}. Similar to \citet[Theorem~1]{Yanikoglu2013} the problem can be made convex by replacing $t_{02} - t_{03}$ in \eqref{eq:Fsubs1} with $u_s^2 + u_s z - u_s \alpha - (v_{01})_s - (v_{04})_s$ and omitting \eqref{eq:Fsubs2} from the problem. Constraint \eqref{eq:Fsubs1} then becomes:
\begin{align}
&& \left( u_s + \frac{z}{2} \right)^2 - u_s \alpha - (v_{01})_s - (v_{04})_s + \rho t_{04} + \sum_{s' \in S} \hat{p}_{s'} \left( \frac{(v_{04})_{s'}^2}{4t_{04}} + (v_{04})_{s'} \right) \leq w, \; \forall s \in S, \label{eq:Fsubsafter}
\end{align}
which is jointly convex in all variables. In order to improve the tractability and accuracy of computing $F(\alpha)$, we cast it as a conic quadratic problem. The only complicating terms are $(v_{04})_s^2 / (4t_{04})$, which can be reformulated using a standard trick. Constraint \eqref{eq:Fsubsafter} is satisfied if and only if there exists auxiliary variables $y_s$ such that the following inequalities are satisfied:
\begin{align*}
&& & \left( u_s + \frac{z}{2} \right)^2 - u_s \alpha - (v_{01})_s - (v_{04})_s + \rho t_{04} + 2 \sum_{s' \in S}\hat{p}_{s'} \left( y_{s'} + (v_{04})_{s'} \right) \leq w, \; \forall s \in S \\
&& & \twonorm{ \begin{pmatrix}(2 v_{04})_s \\ y_s - 4t_{04} \end{pmatrix} } \leq y_s + 4t_{04}, \; \forall s \in S.
\end{align*}
The problem (R-I) can now be solved by determining the root of $F$.

We perform a numerical analysis on $10$ items and a generated data set of $50$ scenarios. In order to incorporate correlations, we first construct a covariance matrix $\vec{A}\vec{A}\transp$, where $\vec{A}$ is a $10 \times 10$ matrix whose entries are uniformly and independently distributed on $[-0.5, 0.5]$. Then, to reflect the idea that a higher risk gives a higher expected return, a vector of expected returns $\vec{\mu}$ is constructed with a linear mapping on the variances of the items. The mapping is constructed such that the item with the smallest variance gets an expected return of 0.01, and the item with the largest variance gets an expected return of 0.20. Finally, the scenarios are drawn, each from a multivariate normal distribution with the constructed mean $\vec{\mu}$ and covariance $\vec{A}\vec{A}\transp$. We solve the model for $\rho=1$, $\hat{p}_s = 0.02$ for all $s$, and $C = 100$ to obtain a robust solution using YALMIP \citep{Lofberg:2010} and MOSEK (Mosek ApS, Denmark) with their default settings. For this value of $\rho$, the probabilities can vary, on average, between 0 and 0.04. Additionally, we solve the same problem for $\rho=0$, i.e., when $p_s = \hat{p}_s$, to obtain a nonrobust solution. We use bisection search on the interval determined by \eqref{eq:lb} and \eqref{eq:ub} and stopped when the interval width was less than $10^{-10}$ . One step in bisection search takes around 2 seconds, of which around 7\% is spent by MOSEK.

The convergence of the bisection method turns out to be adequate. Let $\vec{x_i}$ denote the solution in iteration $i$, and let $\vec{x^*}$ denote the final solution. The initial search interval is $[0.70,20.09]$. The solution $\vec{x_1}$, obtained from solving $F((0.70+20.09)/2)$, is far from optimal: $\infnorm{\vec{x_1} - \vec{x^*}} \approx 2.2$. In each three or four iterations, $\vec{x_i}$ gains one extra digit of accuracy. After 22 iterations, the accuracy has improved to $\infnorm{\vec{x_{22}} - \vec{x^*}} \approx 4.1 \cdot 10^{-7}$. The algorithm terminates after 37 iterations, with no apparent improvement following the 22$^\textrm{nd}$ iteration. Since $\twonorm{\vec{x^*}} \approx 44.9$, the error after 22 iterations is relatively small.

When $p_s = \hat{p}_s$ for all $s$, the mean-variance ratio of the nominal solution (which is 6.34) is indeed lower than that of the robust solution (which is 6.45). For both solutions we determined the worst case $\hat{\vec{p}}$ and the corresponding objective value. The objective of the robust solution (which is 18.62) is slightly better than that of the nominal solution (which is 18.98). This shows that uncertainty may cause a factor three deterioration of the objective value. Relative to this large difference, the difference between the two solutions is small. So, the nominal solution performs quite well for this example. For the worst case probabilities for the robust solution, we have computed the optimal portfolio as if these solutions were known beforehand (perfect hindsight solution). The objective value equals that of the robust solution. So, even though the robust objective could deteriorate substantially, there is no other solution that performs better.

\subsection{Data Envelopment Analysis}\label{sec:dea}
Data Envelopment Analysis (DEA) is a tool to estimate the efficiency of different decision making units (DMUs), based on their inputs and outputs. DEA was originally introduced for not-for-profit companies, e.g., schools where inputs could be number of teacher hours and number of students per class, and outputs could be arithmetic scores and psychological tests of student attitudes, e.g., toward the community \citep{Charnes1978}. The applicability of DEA is not limited to nonprofit organizations. A reference list of more than 4,000 publications on DEA is given by \cite{Emrouznejad2008}.

Let $n_i$ and $n_o$ denote the number of in- and outputs, respectively. The efficiency of a DMU is defined as the largest fraction of weighted outputs divided by weighted inputs, given that the efficiency of the other DMUs is at most 1:
\begin{align*}
\mbox{(DEA)} \quad \max_{\vec{u} \in \R^{n_o}_+,\vec{v} \in \R^{n_i}_+}   \quad \frac{\vec{u} \transp \vec{y_0}}{\vec{v} \transp \vec{x_0}} \quad \mbox{s.t.} \quad \frac{\vec{u} \transp \vec{y_i}}{\vec{v} \transp \vec{x_i}} \leq 1, \; \forall i \in I,
\end{align*}
where $\vec{x_i}$ and $\vec{y_i}$ are the vectors of inputs and outputs of DMU $i$, and $\vec{u}$ and $\vec{v}$ are the non-negative weights.

The inputs and outputs are model parameters that have to be acquired from each DMU and are affected by measurement errors. Especially when a single DMU represents a group of smaller business units or is a pool of all activities in a certain region, and the inputs and outputs are aggregated, errors become practically inevitable. There have been many attempts to incorporate uncertainty in DEA. For an overview, e.g., see \citep{Shokouhi2013}. Since our focus is on RO, we only discuss the three papers that are relevant. The first only considers uncertain outputs \citep{Sadjadi2008}. The second considers jointly uncertain inputs and outputs \citep{Shokouhi2010}. Unfortunately, the robust counterpart in the latter is constructed in an ad-hoc manner that results in a nonconvex formulation, for which it is not clear whether globally optimal solutions were found. The third considers either uncertain inputs or uncertain outputs \citep{Wang2010}. In the last two papers, a simulation study is performed to quantify the improvement offered by the robust solution. For each randomly drawn set of inputs and outputs, they compute the relative efficiencies with the $\vec{u}$ and $\vec{v}$ obtained from the robust solution. However, in our view, when the inputs and outputs are fully known, the relative efficiencies can only be computed by optimizing (DEA) for those known inputs and outputs. Our results are therefore different in two ways. First, we consider both uncertain inputs and uncertain outputs and solve the correct problem. Second, we perform a valid simulation study to verify whether the robust solution is better than the nominal solution.

In this section, we take the data from \cite{Shokouhi2010}. In this data set there are five DMUs, two inputs and two outputs. The in- and outputs are uncertain, but known to reside in given intervals, given in Table \ref{tbl:data:dea}.

In order to get in the minimization framework, the objective of (DEA) is replaced with its reciprocal. The optimal solution of the robust counterpart of (DEA) then corresponds to the reciprocal of the root of:
\begin{align*}
F(\alpha) = \min_{\vec{u} \in \R^{n_o}_+,\vec{v} \in \R^{n_i}_+,w \in \R}  \quad w \quad \mbox{s.t.} \quad & \vec{v} \transp \vec{x_0} - \alpha \vec{u} \transp \vec{y_0} \leq w, \; \forall (\vec{x},\vec{y}) \in \Uset, \\
& \vec{u} \transp \vec{y_i} \leq \vec{v} \transp \vec{x_i}, \; \forall (\vec{x},\vec{y}) \in \Uset, \; \forall i \in I.
\end{align*}
Following \cite{Shokouhi2010}, we take the Bertsimas and Sim uncertainty region:
\begin{align*}
&& \Uset = \{ (\vec{x},\vec{y}) \in \R^{(|I|+1)\times n_i} \times \R^{(|I|+1)\times n_o} : \; & x_{ij} = \bar{x}_{ij} + \zeta^x_{ij} \Delta\hat{x}_{ij}, \quad y_{ij} = \bar{y}_{ij} + \zeta^y_{ij} \Delta\hat{y}_{ij},\\
&&                                   & \infnorm{\mathrm{vec}(\vec{\zeta^x}, \vec{\zeta^y})} \leq 1, \quad \onenorm{\mathrm{vec}(\vec{\zeta^x}, \vec{\zeta^y})} \leq \Gamma \},
\end{align*}
where $\bar{x}_{ij}$ and $\bar{y}_{ij}$ are the midpoints, $\Delta\hat{x}_{ij}$ and $\Delta\hat{y}_{ij}$ are the half-widths of the uncertainty intervals, and the vec operator stacks the columns of the matrix arguments into a single vector. For robust {\it LP}, this set has the property that when $\Gamma$ is integer, it controls the number of uncertain elements that can deviate from their nominal values \citep{BertsimasPoR}. This property also holds for a robust LFP, since $F(\alpha)$ is a robust LP.

The optimal weights $\vec{u}$ and $\vec{v}$ depend on the actual inputs and outputs. One may therefore be inclined to use Adjustable Robust Optimization (ARO) when $\Gamma$ is larger than the dimensions of $\vec{x_i}$ and $\vec{y_i}$ added. Consequently, $\vec{u}$ and $\vec{v}$ are replaced by functions of the uncertain parameters. Unfortunately, even in the simple case of affine decision rules, this is often intractable. In the constraints of $F(\alpha)$, $\vec{u}$ and $\vec{v}$ are multiplied with uncertain parameters, which yields constraints that are quadratic in the uncertain parameters. These can currently only be solved efficiently for ellipsoidal uncertainty sets.

We use bisection search on the interval determined by \eqref{eq:lb} and \eqref{eq:ub} to determine the root of $F(\alpha)$, and stop when the interval width is less than $10^{-4}$ (which turns out to be accurate enough for ranking the DMUs). $F(\alpha)$ is computed using YALMIP and MOSEK with their default settings, and takes a few tenths of a second on a normal desktop computer, where MOSEK accounts for approximately 10\% of that time. The time it takes to compute $F(\alpha)$ turns out to be approximately constant, so independent of the remaining width of the interval and independent of the size of the uncertainty region $\Gamma$. The root of $F(\alpha)$ is determined in a few seconds.

We computed the robust efficiencies of the DMUs for $\Gamma$ ranging between 0 and 4 in steps of 0.1, since each constraint has at most four uncertain parameters. For $\Gamma \leq 0.2$, the list of DMUs ranked from most to least efficient, is 1, 2, 3, 5, 4. For $\Gamma \geq 0.3$, DMUs 3 and 5 switch positions. Hence, DMU 5 is more efficient than DMU 3 when $\Gamma \geq 0.3$. We have tried to verify this claim by running 100 simulations, where in each simulation we uniformly drew inputs and outputs from the uncertainty region, solved (DEA) for each set of inputs and outputs, and ranked the DMUs based on efficiency. In 76 out of 100 simulations, DMU 3 was more efficient than DMU 5. This result advocates against the use of RO in DEA, since it shows that for $\Gamma=0$ (i.e., the nonrobust solution), the ranking is better than for $\Gamma \geq 1$. We have also performed the simulation with more extreme data, by drawing the inputs and outputs only from the endpoints of their uncertainty intervals. This yielded similar results. Other experiments, where we used an ellipsoidal uncertainty region instead of the Bertsimas and Sim uncertainty region, or where we used the nominal objectives (based on $\bar{\vec{x}}$ and $\bar{\vec{y}}$) but kept the uncertain constraints, also yielded similar results.

\begin{table}
	\centering
  \caption{Data set for the DEA example of Section \ref{sec:dea}.}
  \label{tbl:data:dea}
	\begin{tabular}{lllll}
	\toprule
DMU$_i$ & Input 1 & Input 2 & Output 1 & Output 2 \\
   \midrule
1 & [14,15] & [0.06, 0.09] & [157, 161] & [28, 40] \\
2 &  [4,12] & [0.16, 0.35] & [157, 198] & [21, 29] \\
3 & [10,17] & [0.10, 0.70] & [143, 159] & [28, 35] \\
4 & [12,15] & [0.21, 0.48] & [138, 144] & [21, 22] \\
5 & [19,22] & [0.12, 0.19] & [158, 181] & [21, 25] \\
  \bottomrule
	\end{tabular}
\end{table}

\section{Conclusions}
We have shown how RO can be applied to FPs as a method to deal with uncertain data. The method has been tested on three problems. In all three examples, we observe that the nominal solution, which is obtained by solving the deterministic problem, is severely affected by uncertainty. Surprisingly, this also holds for the robust solution, and in none of the examples the robust solution offers a significant improvement; even when comparing worst case performance.

The first question that arises, is why the nominal solution performs so well. We try to answer this question for the mean-variance optimization problem, and note that the explanation for the multi-item newsvendor problem is similar. For a given solution, the worst case for the mean is when the probability vector is a unit vector, that assigns unit weight to the scenario with lowest return. For the variance, the worst case is when the scenarios with the lowest and highest return each occur with probability 0.5. For a robust solution w.r.t. the mean value, the scenario with lowest return should be optimized, whereas for the variance, the returns in the scenarios with lowest and highest return should be close to each other. The nominal solution simultaneously maximizes the expected value and minimizes the variance. While not identical to the robust objective, it contains some aspects of it. For example, the mean is a weighted sum that contains the return for the scenario with lowest return.

The second question that arises, is why there is a realization of the uncertain parameters in Sections \ref{sec:mine} and \ref{sec:mvo}, for which no solution can outperform the robust solution; even if the former is optimized as if the realization of the uncertain parameters are known beforehand. This turns out to be due to Sion's minimax theorem \citep{Sion1958}. The assumptions (a)-(g) ensure that $f(\vec{a_0}, \vec{x}) / g(\vec{a_0}, \vec{x})$ is quasi-convex in $\vec{x}$ (for fixed $\vec{a_0}$), quasi-concave in $\vec{a_0}$ (for fixed $\vec{x}$), and continuous, that the uncertainty set is compact and convex, and that the feasible set for $\vec{x}$, say $X$, is convex. Therefore, by Sion's minimax theorem, $\max_{\vec{a_0} \in \Uset} \min_{\vec{x} \in X} f(\vec{a_0}, \vec{x}) / g(\vec{a_0}, \vec{x}) = \min_{\vec{x} \in X} \max_{\vec{a_0} \in \Uset} f(\vec{a_0}, \vec{x}) / g(\vec{a_0}, \vec{x})$. This no longer holds when there is uncertainty in the constraints, since the feasible region $X$ changes when the values for the uncertain parameters are known.

So, the robust solution is good in the sense that it cannot be improved in the worst case, even if the values of the uncertain parameters are known beforehand. On the other hand, the nominal solution performs well, at least in the examples studied. It shall be interesting to see the difference in real-life examples, especially with uncertainty in the constraints.

\subsubsection*{Acknowledgments}
We thank D. den Hertog (Tilburg University) for many useful ideas and comments, D. Iancu (Stanford University) for showing the formulation (RP-FP) used in Theorem \ref{thm:dual}, and the referee and editor for their valuable comments.

\bibliographystyle{abbrvnatnew}
\bibliography{library}

\begin{thebibliography}{34}
\providecommand{\natexlab}[1]{#1}
\providecommand{\url}[1]{\texttt{#1}}
\expandafter\ifx\csname urlstyle\endcsname\relax
  \providecommand{\doi}[1]{doi: #1}\else
  \providecommand{\doi}{doi: \begingroup \urlstyle{rm}\Url}\fi

\bibitem[Barros et~al.(1996)Barros, Frenk, Schaible, and Zhang]{Barros1996}
A.~I. Barros, J.~B.~G. Frenk, S.~Schaible, and S.~Zhang.
\newblock A new algorithm for generalized fractional programs.
\newblock \emph{Mathematical Programming},
  \href{http://dx.doi.org/10.1007/BF02592087}{72\penalty0 (2):\penalty0
  147--175}, 1996.

\bibitem[Beck and Ben-Tal(2009)]{Beck2009}
A.~Beck and A.~Ben-Tal.
\newblock {Duality in robust optimization: primal worst equals dual best}.
\newblock \emph{Operations Research Letters},
  \href{http://dx.doi.org/10.1016/j.orl.2008.09.010}{37\penalty0 (1):\penalty0
  1--6}, 2009.

\bibitem[Ben-Tal and den Hertog(2014)]{BenTal2014hidden}
A.~Ben-Tal and D.~den Hertog.
\newblock Hidden conic quadratic representation of some nonconvex quadratic
  optimization problems.
\newblock \emph{Mathematical Programming},
  \href{http://dx.doi.org/10.1007/s10107-013-0710-8}{143\penalty0
  (1-2):\penalty0 1--29}, 2014.

\bibitem[Ben-Tal et~al.(2009)Ben-Tal, {El Ghaoui}, and Nemirovski]{BenTal}
A.~Ben-Tal, L.~{El Ghaoui}, and A.~Nemirovski.
\newblock \emph{Robust Optimization}.
\newblock Princeton Series in Applied Mathematics. Princeton University Press,
  2009.

\bibitem[Ben-Tal et~al.(2013)Ben-Tal, den Hertog, de~Waegenaere, Melenberg, and
  Rennen]{BTHWMR2011}
A.~Ben-Tal, D.~den Hertog, A.~M.~B. de~Waegenaere, B.~Melenberg, and G.~Rennen.
\newblock Robust solutions of optimization problems affected by uncertain
  probabilities.
\newblock \emph{Management Science},
  \href{http://dx.doi.org/10.1287/mnsc.1120.1641}{59\penalty0 (2):\penalty0
  341--357}, 2013.

\bibitem[Ben-Tal et~al.(2015)Ben-Tal, den Hertog, and \noop{J.-Ph.}
  Vial]{BenTal2011}
A.~Ben-Tal, D.~den Hertog, and \noop{J.-Ph.} Vial.
\newblock Deriving robust counterparts of nonlinear uncertain inequalities.
\newblock \emph{Mathematical Programming},
  \href{http://dx.doi.org/10.1007/s10107-014-0750-8}{149\penalty0
  (1--2):\penalty0 265--299}, 2015.

\bibitem[Bertsimas and Sim(2004)]{BertsimasPoR}
D.~Bertsimas and M.~Sim.
\newblock The price of robustness.
\newblock \emph{Operations Research},
  \href{http://dx.doi.org/10.1287/opre.1030.0065}{52\penalty0 (1):\penalty0
  35--53}, 2004.

\bibitem[Bertsimas et~al.(2011{\natexlab{a}})Bertsimas, Brown, and
  Caramanis]{Bertsimas2011}
D.~Bertsimas, D.~B. Brown, and C.~Caramanis.
\newblock Theory and applications of robust optimization.
\newblock \emph{SIAM Review},
  \href{http://dx.doi.org/10.1137/080734510}{53\penalty0 (3):\penalty0
  464--501}, 2011{\natexlab{a}}.

\bibitem[Bertsimas et~al.(2011{\natexlab{b}})Bertsimas, Iancu, and
  Parrilo]{BertsimasPolynomials}
D.~Bertsimas, D.~A. Iancu, and P.~A. Parrilo.
\newblock A hierarchy of near-optimal policies for multi-stage adaptive
  optimization.
\newblock \emph{IEEE Transactions on Automatic Control},
  \href{http://dx.doi.org/10.1287/ijoc.1100.0419}{56\penalty0 (12):\penalty0
  2809--2824}, 2011{\natexlab{b}}.

\bibitem[Charnes and Cooper(1962)]{CharnesCooper}
A.~Charnes and W.~W. Cooper.
\newblock Programming with linear fractional functionals.
\newblock \emph{Naval Research Logistics Quarterly},
  \href{http://dx.doi.org/10.1002/nav.3800090303}{9\penalty0 (3--4):\penalty0
  181--186}, 1962.

\bibitem[Charnes et~al.(1978)Charnes, Cooper, and Rhodes]{Charnes1978}
A.~Charnes, W.~W. Cooper, and E.~Rhodes.
\newblock Measuring the efficiency of decision making units.
\newblock \emph{European Journal of Operational Research},
  \href{http://dx.doi.org/10.1016/0377-2217(78)90138-8}{2\penalty0
  (6):\penalty0 429--444}, 1978.

\bibitem[Chen et~al.(2009)Chen, Schaible, and Sheu]{ChenFP2009}
H.~J. Chen, S.~Schaible, and R.~L. Sheu.
\newblock Generic algorithm for generalized fractional programming.
\newblock \emph{Journal of Optimization Theory and Applications},
  \href{http://dx.doi.org/10.1007/s10957-008-9499-7}{141\penalty0 (1):\penalty0
  93--105}, 2009.

\bibitem[Crouzeix and Ferland(1991)]{crouzeix1991}
J.-P. Crouzeix and J.~A. Ferland.
\newblock Algorithms for generalized fractional programming.
\newblock \emph{Mathematical Programming},
  \href{http://dx.doi.org/10.1007/BF01582887}{52\penalty0 (1--3):\penalty0
  191--207}, 1991.

\bibitem[Crouzeix et~al.(1985)Crouzeix, Ferland, and Schaible]{crouzeix1985}
J.-P. Crouzeix, J.~A. Ferland, and S.~Schaible.
\newblock An algorithm for generalized fractional programs.
\newblock \emph{Journal of Optimization Theory and Applications},
  \href{http://dx.doi.org/10.1007/BF00941314}{47\penalty0 (1):\penalty0
  35--49}, 1985.

\bibitem[Crouzeix et~al.(1986)Crouzeix, Ferland, and Schaible]{crouzeix1986}
J.-P. Crouzeix, J.~A. Ferland, and S.~Schaible.
\newblock A note on an algorithm for generalized fractional programs.
\newblock \emph{Journal of Optimization Theory and Applications},
  \href{http://dx.doi.org/10.1007/BF00938484}{50\penalty0 (1):\penalty0
  183--187}, 1986.

\bibitem[Dinkelbach(1967)]{Dinkelbach1967}
W.~Dinkelbach.
\newblock On nonlinear fractional programming.
\newblock \emph{Management Science},
  \href{http://dx.doi.org/10.1287/mnsc.13.7.492}{13\penalty0 (7):\penalty0
  492--498}, 1967.

\bibitem[Emrouznejad et~al.(2008)Emrouznejad, Parker, and
  Tavares]{Emrouznejad2008}
A.~Emrouznejad, B.~R. Parker, and G.~Tavares.
\newblock Evaluation of research in efficiency and productivity: A survey and
  analysis of the first 30 years of scholarly literature in {DEA}.
\newblock \emph{Socio-Economic Planning Sciences},
  \href{http://dx.doi.org/10.1016/j.seps.2007.07.002}{42\penalty0 (3):\penalty0
  151--157}, 2008.

\bibitem[Gorissen et~al.(2014)Gorissen, Ben-Tal, Blanc, and den
  Hertog]{Gorissen2013convunc}
B.~L. Gorissen, A.~Ben-Tal, J.~P.~C. Blanc, and D.~den Hertog.
\newblock Deriving robust and globalized robust solutions of uncertain linear
  programs with general convex uncertainty sets.
\newblock \emph{Operations Research},
  \href{http://dx.doi.org/10.1287/opre.2014.1265}{62\penalty0 (3):\penalty0
  672--679}, 2014.

\bibitem[Jeyakumar et~al.(2013)Jeyakumar, Li, and
  Srisatkunarajah]{Jeyakumar2013}
V.~Jeyakumar, G.~Y. Li, and S.~Srisatkunarajah.
\newblock Strong duality for robust minimax fractional programming problems.
\newblock \emph{European Journal of Operational Research},
  \href{http://dx.doi.org/10.1016/j.ejor.2013.02.015}{228\penalty0
  (2):\penalty0 331--336}, 2013.

\bibitem[Kaul et~al.(1986)Kaul, Kaur, and Lyall]{Kaul1986}
R.~N. Kaul, S.~Kaur, and V.~Lyall.
\newblock Duality in inexact fractional programming with set-inclusive
  constraints.
\newblock \emph{Journal of Optimization Theory and Applications},
  \href{http://dx.doi.org/10.1007/BF00939274}{50\penalty0 (2):\penalty0
  279--288}, 1986.

\bibitem[Lin and Sheu(2005)]{Lin2005}
J.~Lin and R.~Sheu.
\newblock Modified {D}inkelbach-type algorithm for generalized fractional
  programs with infinitely many ratios.
\newblock \emph{Journal of Optimization Theory and Applications},
  \href{http://dx.doi.org/10.1007/s10957-005-4717-z}{126\penalty0 (2):\penalty0
  323--343}, 2005.

\bibitem[L\"{o}fberg(2012)]{Lofberg:2010}
J.~L\"{o}fberg.
\newblock Automatic robust convex programming.
\newblock \emph{Optimization Methods and Software},
  \href{http://dx.doi.org/10.1080/10556788.2010.517532}{27\penalty0
  (1):\penalty0 115--129}, 2012.

\bibitem[Markowitz(1952)]{Markowitz1952}
H.~M. Markowitz.
\newblock Portfolio selection.
\newblock \emph{Journal of Finance},
  \href{http://dx.doi.org/10.2307/2975974}{7\penalty0 (1):\penalty0 77--91},
  1952.

\bibitem[Sadjadi and Omrani(2008)]{Sadjadi2008}
S.~J. Sadjadi and H.~Omrani.
\newblock Data envelopment analysis with uncertain data: an application for
  {I}ranian electricity distribution companies.
\newblock \emph{Energy Policy},
  \href{http://dx.doi.org/10.1016/j.enpol.2008.08.004}{36\penalty0
  (11):\penalty0 4247--4254}, 2008.

\bibitem[Schaible(1974)]{Schaible1974}
S.~Schaible.
\newblock Parameter-free convex equivalent and dual programs of fractional
  programming problems.
\newblock \emph{Zeitschrift für Operations Research},
  \href{http://dx.doi.org/10.1007/BF02026600}{18\penalty0 (5):\penalty0
  187--196}, 1974.

\bibitem[Schaible(1976)]{Schaible1976}
S.~Schaible.
\newblock Fractional programming. {II}, {O}n {D}inkelbach's algorithm.
\newblock \emph{Management Science},
  \href{http://dx.doi.org/10.1287/mnsc.22.8.868}{22\penalty0 (8):\penalty0
  868--873}, 1976.

\bibitem[Schaible(1982)]{Schaible1982}
S.~Schaible.
\newblock Bibliography in fractional programming.
\newblock \emph{Zeitschrift für Operations Research},
  \href{http://dx.doi.org/10.1007/BF01917115}{26\penalty0 (1):\penalty0
  211--241}, 1982.

\bibitem[Schaible and Ibaraki(1983)]{Schaible1983}
S.~Schaible and T.~Ibaraki.
\newblock Fractional programming.
\newblock \emph{European Journal of Operational Research},
  \href{http://dx.doi.org/10.1016/0377-2217(83)90153-4}{12\penalty0
  (4):\penalty0 325--338}, 1983.

\bibitem[Shokouhi et~al.(2010)Shokouhi, Hatami-Marbini, Tavana, and
  Saati]{Shokouhi2010}
A.~H. Shokouhi, A.~Hatami-Marbini, M.~Tavana, and S.~Saati.
\newblock A robust optimization approach for imprecise data envelopment
  analysis.
\newblock \emph{Computers \& Industrial Engineering},
  \href{http://dx.doi.org/10.1016/j.cie.2010.05.011}{59\penalty0 (3):\penalty0
  387--397}, 2010.

\bibitem[Shokouhi et~al.(2014)Shokouhi, Shahriari, Agrell, and
  Hatami-Marbini]{Shokouhi2013}
A.~H. Shokouhi, H.~Shahriari, P.~Agrell, and A.~Hatami-Marbini.
\newblock Consistent and robust ranking in imprecise data envelopment analysis
  under perturbations of random subsets of data.
\newblock \emph{OR Spectrum},
  \href{http://dx.doi.org/10.1007/s00291-013-0336-5}{36\penalty0 (1):\penalty0
  133--160}, 2014.

\bibitem[Sion(1958)]{Sion1958}
M.~Sion.
\newblock On general minimax theorems.
\newblock \emph{Pacific Jouronal of Mathematics},
  \href{http://projecteuclid.org/euclid.pjm/1103040253}{8\penalty0
  (1):\penalty0 171--176}, 1958.

\bibitem[Stancu-Minasian(2013)]{Stancu2013}
I.~M. Stancu-Minasian.
\newblock A seventh bibliography of fractional programming.
\newblock \emph{Advanced Modeling and Optimization},
  \href{http://camo.ici.ro/journal/vol15/v15b13.pdf}{15\penalty0 (2):\penalty0
  309--386}, 2013.

\bibitem[Wang and Wei(2010)]{Wang2010}
K.~Wang and F.~Wei.
\newblock Robust data envelopment analysis based {MCDM} with the consideration
  of uncertain data.
\newblock \emph{Journal of Systems Engineering and Electronics},
  \href{http://dx.doi.org/10.3969/j.issn.1004-4132.2010.06.009}{21\penalty0
  (6):\penalty0 981--989}, 2010.

\bibitem[Yan{\i}ko\u{g}lu et~al.(2013)Yan{\i}ko\u{g}lu, den Hertog, and
  Kleijnen]{Yanikoglu2013}
{\.{I}}.~Yan{\i}ko\u{g}lu, D.~den Hertog, and J.~P.~C. Kleijnen.
\newblock Adjustable robust parameter design with unknown distributions.
\newblock \emph{CentER Discussion Paper},
  \href{http://arno.uvt.nl/show.cgi?fid=129316}{2013\penalty0 (022)}, 2013.

\end{thebibliography}

\appendix
\section{The Importance of Convexity Conditions}\label{sec:lin2005}
We provide a short example to stress the importance of the convexity/concavity conditions on $f$ and $g$. The second numerical example by \cite{Lin2005} is:
\begin{align*}
\min_{\vec{x} \in \R^n} \max_{a \in [0,1]} \quad \frac{ a^2 x_1 x_2 + x_1^{2a} + a x_3^3  }{  5(a-1)^2 x_1^4 + 2x_2^2 + 4ax_3} \quad \mbox{s.t.} \quad 0.5 \leq x_i \leq 5, \quad i=1,2,3.
\end{align*}
This problem does not satisfy the convexity/concavity conditions from Section \ref{sec:solving-rfp}. Lin and Sheu claim that $\vec{x} = (0.5 \quad 1.5 \quad 0.5)$ and $a=0$ is optimal with a value of 0.21 (reported as $-0.21$), but $\vec{x} = (0.5 \quad 5 \quad 0.5)$ and $a=1$ is a better solution (maybe still not optimal) since the corresponding value is 0.06.

\section{On the Result by Kaul et al. (1986)}\label{sec:kaul86}
This appendix shows a mistake in the paper by \cite{Kaul1986}. Essentially, they formulate the dual of:
\begin{align*}
\min_{\alpha \in \R_+,\vec{x} \in \R^n_+} \; \alpha  \quad \mbox{s.t.} \quad \frac{b_0 + \vec{b} \transp \vec{x}}{c_0 + \vec{c} \transp \vec{x}} \leq \alpha, \; \forall (b_0,\vec{b}) \in \Uset_{1} \times \Uset_{2}, \; \forall (c_0,\vec{c}) \in \Uset_{3} \times \Uset_{4}, \quad \vec{A} \vec{x} \leq \vec{d}.
\end{align*}
Note that $\vec{x}$ is non-negative. In their Lemma 2.1, they claim that the worst case $(c_0,\vec{c})$ does not depend on $\vec{x}$, and is given by $c_0^* = \min_{c_0 \in \Uset_{3}} \{ c_0 \}$ and $\vec{c^*}$, with components $c_i^* = \min_{\vec{c} \in \Uset_{4}} \{ c_i \}$. This implicitly assumes that $\vec{c^*}$ is a member of $\Uset_{4}$, which is not always true. The mistake becomes clear in their numerical example, where they use $\vec{c^*} = [4; 2]$, which is not in the uncertainty set. Consequently, the proposed approach gives the wrong dual problem and a suboptimal solution. Our results in Section \ref{sec:singlesolverrun} can provide the correct dual problem under milder conditions on the uncertainty sets.

\end{document}